\documentclass[a4paper]{article}
\usepackage{spconf,amsmath,cite}
\usepackage{amssymb}
\usepackage{euscript}
\usepackage{comment}
\usepackage{theorem}
\usepackage{charter}
\usepackage{mathrsfs} 
\usepackage{euscript}
\usepackage{graphicx}
\usepackage{xcolor}
\usepackage{enumitem}
\usepackage{url}
\definecolor{labelkey}{rgb}{0,0.08,0.45}
\definecolor{refkey}{rgb}{0,0.6,0.0}
\definecolor{Brown}{rgb}{0.45,0.0,0.05}
\definecolor{dgreen}{rgb}{0.00,0.49,0.00}
\definecolor{dblue}{rgb}{0,0.08,0.75}
\RequirePackage[colorlinks,hyperindex]{hyperref}
\hypersetup{linktocpage=true,citecolor=dblue,linkcolor=dgreen}

\renewcommand{\leq}{\ensuremath{\leqslant}}

\renewcommand{\le}{\ensuremath{\leqslant}}

\newcommand{\minimize}[2]{\ensuremath{\underset{\substack{{#1}}}%
{\text{\rm minimize}}\;\;#2 }}
\newcommand{\EC}[2]{{\mathsf E}(#1\! \mid\! #2)}

\newcommand{\scal}[2]{{\left\langle{{#1}\mid{#2}}\right\rangle}}

\newcommand{\menge}[2]{\big\{{#1}~\big |~{#2}\big\}} 
\newcommand{\Menge}[2]{\left\{{#1}~\Big|~{#2}\right\}} 

\newcommand{\GGG}{{\ensuremath{\boldsymbol{\mathsf G}}}}

\newcommand{\HH}{\ensuremath{{\mathsf H}}}
\newcommand{\GG}{\ensuremath{{\mathsf G}}}
\newcommand{\FF}{\ensuremath{{\mathcal F}}}
\newcommand{\XX}{\ensuremath{\EuScript{X}}}

\newcommand{\XXX}{\ensuremath{\boldsymbol{\EuScript{X}}}}

\newcommand{\Sum}{\ensuremath{\displaystyle\sum}}

\newcommand{\RR}{\ensuremath{\mathbb{R}}}
\newcommand{\RP}{\ensuremath{\left[0,+\infty\right[}}

\newcommand{\RPP}{\ensuremath{\left]0,+\infty\right[}}

\newcommand{\RX}{\ensuremath{\left]-\infty,+\infty\right]}}

\newcommand{\EE}{\ensuremath{\mathsf E}}
\newcommand{\PP}{\ensuremath{\mathsf P}}
\newcommand{\as}{\ensuremath{\text{\rm $\PP$-a.s.}}}

\newcommand{\NN}{\ensuremath{\mathbb N}}

\newcommand{\pinf}{\ensuremath{{+\infty}}}

\newcommand{\prox}{\ensuremath{\text{\rm prox}}}

\newcommand{\infconv}{\ensuremath{\mbox{\small$\,\square\,$}}}

\newcommand{\rzeroun}{\ensuremath{\left]0,1\right]}}   
\newcommand{\argmind}[2]{\ensuremath{\underset{\substack{{#1}}}%
{\text{argmin}}\;\left(#2\right)}}

\newtheorem{theorem}{Theorem}[section]

\newtheorem{proposition}[theorem]{Proposition}
\theoremstyle{plain}{\theorembodyfont{\rmfamily}%
}
\theoremstyle{plain}{\theorembodyfont{\rmfamily}%
}
\theoremstyle{plain}{\theorembodyfont{\rmfamily}%
\newtheorem{algorithm}[theorem]{Algorithm}}
\theoremstyle{plain}{\theorembodyfont{\rmfamily}%
\newtheorem{problem}[theorem]{Problem}}
\theoremstyle{plain}{\theorembodyfont{\rmfamily}%
}
\theoremstyle{plain}{\theorembodyfont{\rmfamily}%
}
\theoremstyle{plain}{\theorembodyfont{\rmfamily}%
}
\theoremstyle{plain}{\theorembodyfont{\rmfamily}%
}

\title{STOCHASTIC FORWARD-BACKWARD AND PRIMAL-DUAL APPROXIMATION
ALGORITHMS WITH APPLICATION TO ONLINE IMAGE RESTORATION}
\name{Patrick L. Combettes$^{1}$ and 
Jean-Christophe Pesquet$^2$\thanks{This work was supported by the
CNRS Imag'In project under grant 2015 OPTIMISME.}}
\address{
\normalsize $\!^1$Sorbonne Universit\'es -- UPMC Univ. Paris 06,
UMR 7598, Laboratoire Jacques-Louis Lions,\\
\normalsize Paris, France,
\normalsize plc@ljll.math.upmc.fr\\
\normalsize $\!^2$Universit\'e Paris-Est, Laboratoire
d'Informatique Gaspard Monge -- CNRS UMR 8049,\\
\normalsize Champs-sur-Marne, France,
\normalsize jean-christophe.pesquet@univ-paris-est.fr
}

\begin{document}

\maketitle
\begin{abstract}
Stochastic approximation techniques have been used in various
contexts in data science. We propose a stochastic version of the
forward-backward algorithm for minimizing the sum of two
convex functions, one of which is not necessarily smooth. Our
framework can handle stochastic approximations of the
gradient of the smooth function and allows for stochastic errors in
the evaluation of the proximity operator of the nonsmooth function.
The almost sure convergence of the iterates generated
by the algorithm to a minimizer 
is established under relatively mild assumptions. We also
propose a stochastic version of a popular primal-dual proximal
splitting algorithm, establish its convergence, and apply it to
an online image restoration problem.
\end{abstract}
\begin{keywords}
convex optimization, nonsmooth optimization, primal-dual algorithm,
stochastic algorithm, parallel algorithm, proximity operator,
recovery, image restoration.
\end{keywords}
\section{Introduction}
\label{sec:intro}


A large array of optimization problems arising in signal processing 
involve functions belonging to $\Gamma_0(\HH)$, the class of 
proper lower semicontinuous convex function from $\HH$ to $\RX$, 
where $\HH$ is a finite-dimensional real Hilbert space with 
norm $\|\cdot\|$. In particular, the following formulation has
proven quite flexible and far reaching \cite{Smms05}.

\begin{problem}
\label{prob:1}
Let $\mathsf{f}\in \Gamma_0(\HH)$, let $\vartheta\in\RPP$, and let 
$\mathsf{g}\colon\HH\to\RR$ be a differentiable convex function 
such that $\nabla\mathsf{g}$ is $\vartheta^{-1}$-Lipschitz 
continuous on $\HH$. The goal is to
\begin{equation}
\label{e:2015-01-20p}
\minimize{\mathsf{x}\in\HH}
{\mathsf{f}(\mathsf{x})+\mathsf{g}(\mathsf{x})},
\end{equation}
under the assumption that the set
$\mathsf{F}$ of minimizers of $\mathsf{f}+\mathsf{g}$ is nonempty.
\end{problem}

A standard method to solve Problem~\ref{prob:1} is the 
forward-backward algorithm 
\cite{Smms05,Chau07,Siop07,Banf11,Yama15}, 
which constructs a sequence $(\mathsf{x}_n)_{n\in\NN}$ in 
$\HH$ via the recursion
\begin{equation}
\label{e:2015-02-14a}
(\forall n\in\NN) \quad \mathsf{x}_{n+1}=
\prox_{\gamma_n\mathsf{f}}\big(\mathsf{x}_n-\gamma_n
\mathsf{\nabla}\mathsf{g}(\mathsf{x}_n)\big),
\end{equation}
where $\gamma_n\in\left]0,2\vartheta\right[$ 
and $\prox_{\gamma_n\mathsf{f}}$ 
is the proximity operator of function $\gamma_n \mathsf{f}$, i.e.,
\cite{Livre1}
\begin{equation}
\prox_{\gamma_n\mathsf{f}}\colon
\mathsf{x}\to
\argmind{\mathsf{y}\in\HH}{\mathsf{f}(\mathsf{y})+
\frac{1}{2\gamma_n}\|\mathsf{x}-\mathsf{y}\|^2}.
\end{equation}
In practice, it may happen that, at each iteration $n$,
$\mathsf{\nabla}\mathsf{g}(\mathsf{x}_n)$ is not known exactly and
is available only through some stochastic approximation $u_n$,
while only a deterministic approximation $\mathsf{f}_n$ to
$\mathsf{f}$ is known; see, e.g., \cite{Pey16}. To solve 
\eqref{e:2015-01-20p} in such uncertain environments, we propose to 
investigate the following stochastic version of 
\eqref{e:2015-02-14a}.
In this algorithm, at iteration $n$, $a_n$ stands for a 
stochastic error term modeling inexact implementations of the 
proximity operator of $\gamma_n \mathsf{f}_n$, 
$(\Omega,\FF,\PP)$ is the underlying probability space, 
and $L^2(\Omega,\FF,\PP;\HH)$ denotes the space
of $\HH$-valued random variable $x$ such that $\EE\|x\|^2<\pinf$.
Our algorithmic model is the following.

\begin{algorithm}
\label{algo:1}
Let $x_0$, $(u_n)_{n\in\NN}$, and $(a_n)_{n\in\NN}$ be random 
variables in $L^2(\Omega,\FF,\PP;\HH)$, let 
$(\lambda_n)_{n\in\NN}$ be a sequence in $\rzeroun$, and let
$(\gamma_n)_{n\in\NN}$ be a sequence in 
$\left]0,2\vartheta\right[$, and let 
$(\mathsf{f}_n)_{n\in\NN}$ be a sequence of functions in
$\Gamma_0(\HH)$. For every $n\in\NN$, set
\begin{equation}
\label{e:FB}
x_{n+1} =
x_n+\lambda_n\big(\prox_{\gamma_n \mathsf{f}_n}
(x_n-\gamma_nu_n)+a_n-x_n\big).
\end{equation}
\end{algorithm}

The first instances of the stochastic iteration \eqref{e:FB} can be 
traced back to \cite{Robi51} in the context of the gradient descent
method, i.e., when $\mathsf{f}_n \equiv \mathsf{f} =\mathsf{0}$. 
Stochastic approximations in the gradient method 
were then investigated in the Russian literature of the late 
1960s and early 1970s 
\cite{Ermo67,Guse71,Shor85}. 
Stochastic gradient methods have also been used extensively in 
adaptive signal processing, in control, and in machine learning,
(e.g., in \cite{Bach11,Kush03,Widr03}). More generally, proximal 
stochastic gradient methods have been applied to various problems; 
see for instance 
\cite{Atch14,Duch09,Kon14,Ros14b,SS2013,Xiao14,Yamagishi11}. 

The first objective of the present work is to provide a thorough
convergence analysis of the stochastic forward-backward 
algorithm  described in Algorithm~\ref{algo:1}.
In particular, our results do 
not require that the proximal parameter sequence 
$(\gamma_n)_{n\in\NN}$ be vanishing.
A second goal of our paper is to show that the extension of
Algorithm~\ref{algo:1} for solving monotone inclusion problems
allows us to derive a stochastic version of a recent primal-dual
algorithm \cite{Bang13} (see also \cite{Opti14,Cond13}). 
Note that our 
algorithm is different from the random block-coordinate approaches 
developed in \cite{Bian14,Repe15}, and that it is more in the 
spirit of the adaptive method of \cite{Ono13}.

The organization of the paper is as follows. Section~\ref{sec:FB}
contains our main result on the convergence of the iterates of
Algorithm~\ref{algo:1}. Section~\ref{sec:PD} presents a 
stochastic primal-dual approach for solving composite convex
optimization problems. Section~\ref{se:appli} illustrates the
benefits of this algorithm in signal restoration problems with 
stochastic degradation operators. Concluding
remarks appear in Section~\ref{se:conclu}.

\section{A stochastic forward-backward algorithm}
\label{sec:FB}
Throughout, given a sequence $(x_n)_{n\in\NN}$ of $\HH$-valued 
random variables, the smallest $\sigma$-algebra generated by 
$x_0,\ldots,x_n$ is denoted by $\sigma(x_0,\ldots,x_n)$, and
we denote by 
$\mathscr{X}=(\XX_n)_{n\in\NN}$ a sequence
of sigma-algebras such that
\begin{equation}
\label{e:2013-11-14}
(\forall n\in\NN)\quad\XX_n\subset\FF\quad\text{and}\quad
\sigma(x_0,\ldots,x_n)\subset\XX_n\subset\XX_{n+1}.
\end{equation}
Furthermore, $\ell_+(\mathscr{X})$ designates the set of sequences 
of $\RP$-valued random variables $(\xi_n)_{n\in\NN}$ such that,
for every $n\in\NN$, $\xi_n$ is $\XX_n$-measurable, and we define
\begin{equation}
\label{e:2013-11-13}
\ell_+^{1/2}(\mathscr{X})=
\Menge{(\xi_n)_{n\in\NN}\in\ell_+(\mathscr{X})}
{\sum_{n\in\NN}\xi_n^{1/2}<\pinf\;\as},
\end{equation}
and
\begin{equation}
\label{e:2013-11-12}
\ell_+^\infty({\mathscr{X}})=
\Menge{(\xi_n)_{n\in\NN}\in\ell_+(\mathscr{X})}
{\sup_{n\in\NN}\xi_n<\pinf\; \as}.
\end{equation}
We now state our main convergence result.

\begin{theorem}
\label{t:2}
Consider the setting of Problem~\ref{prob:1}, let
$(\tau_n)_{n\in\NN}$ be a sequence in $\RP$, 
let $(x_n)_{n\in\NN}$ be a sequence generated by 
Algorithm~\ref{algo:1}, and
let $\mathscr{X}=(\XX_n)_{n\in\NN}$ be a sequence of 
sub-sigma-algebras satisfying \eqref{e:2013-11-14}. 
Suppose that the following are satisfied:
\begin{enumerate}[label=\rm(\alph*)]
\item 
\label{a:t2i}
$\sum_{n\in\NN}\lambda_n\sqrt{\EC{\|a_n\|^2}{\XX_n}}<\pinf$.
\item 
\label{a:t2ii}
$\sum_{n\in\NN}\sqrt{\lambda_n}
\|\EC{u_n}{\XX_n}-\nabla\mathsf{g}(x_n)\|<\pinf$.
\item 
\label{a:t2iii} 
For every $\mathsf{z}\in\mathsf{F}$, there exists 
$(\zeta_n(\mathsf{z}))_{n\in\NN}\in
\ell^\infty_+({\mathscr{X}})$ such that
$\big(\lambda_n\zeta_n(\mathsf{z})\big)_{n\in\NN}\in
\ell_+^{1/2}({\mathscr{X}})$ and
\begin{multline}
\label{e:boundtauetazeta}
(\forall n\in\NN)\quad 
\EC{\|u_n-\EC{u_n}{\XX_n}\|^2}{\XX_n}\\
\leq\tau_n\|\nabla\mathsf{g}(x_n)
-\nabla\mathsf{g}(\mathsf{z})\|^2+\zeta_n(\mathsf{z}).
\end{multline}
\item 
\label{p:23vi}  
There exist sequences $(\alpha_n)_{n\in \NN}$ and 
$(\beta_n)_{n\in\NN}$ in $\RP$ such that 
$\sum_{n\in\NN}\sqrt{\lambda_n} \alpha_n < \pinf$, 
$\sum_{n\in\NN}\lambda_n\beta_n$ $<\pinf$, and
\begin{multline}
\label{e:approxprox}
(\forall n\in\NN)(\forall\mathsf{x}\in\HH)\\
\|\prox_{\gamma_n\mathsf{f}_n}\mathsf{x}-
\prox_{\gamma_n \mathsf{f}}\mathsf{x}\| 
\leq\alpha_n \|\mathsf{x}\|+\beta_n.
\end{multline}
\item 
\label{a:t2iv} 
$\inf_{n\in\NN}\gamma_n>0$, $\sup_{n\in\NN}\tau_n<\pinf$, 
and\linebreak
$\sup_{n\in\NN}(1+\tau_n)\gamma_n<2\vartheta$.
\item 
\label{a:t2v} 
Either $\inf_{n\in\NN}\lambda_n>0$ or
$\big[\,\gamma_n\equiv\gamma$, $\sum_{n\in\NN}\tau_n<\pinf$, and
$\sum_{n\in\NN}\lambda_n=\pinf\,\big]$.
\end{enumerate}
Then the following hold for every $\mathsf{z} \in \mathsf{F}$ and for some $\mathsf{F}$-valued random 
variable $x$:
\begin{enumerate}
\item
\label{t:2i}
$\sum_{n\in\NN}\lambda_n\|\nabla\mathsf{g}(x_n)
-\nabla \mathsf{g}(\mathsf{z})\|^2
<\pinf\;\as$
\item
\label{t:2ii}
$\sum_{n\in\NN}\lambda_n\|x_n-\gamma_n\nabla\mathsf{g}(x_n)-
\prox_{\gamma_n\mathsf{f}}\big(x_n-\gamma_n\nabla\mathsf{g}(x_n)\big)
+\gamma_n\nabla\mathsf{g}(\mathsf{z})\|^2<\pinf\;\as$
\item
\label{t:2iii}
$(x_n)_{n\in\NN}$ converges almost surely to $x$.
\end{enumerate}
\end{theorem}

In the deterministic case, Theorem~\ref{t:2}\ref{t:2iii} can be
found in \cite[Corollary~6.5]{Opti04}.
The proof the above stochastic version is based on the
theoretical tools of \cite{Siop15} (see \cite{Pafa1} for technical
details and extensions to infinite-dimensional Hilbert spaces).

It should be noted that the existing works which are the most
closely related to ours do not allow any approximation of the
function $\mathsf{f}$ and make some additional restrictive
assumptions. For example, in \cite[Corollary~8]{Atch14} and
\cite{Ros14a},  $(\gamma_n)_{n\in\NN}$ is a decreasing sequence. In
\cite[Corollary~8]{Atch14}, \cite{Ros14a}, and \cite{Rosa15}, no
error term is allowed in the numerical evaluations of the proximity
operators ($a_n \equiv 0$). 
In addition, in the former work, it is assumed that
$(x_n)_{n\in\NN}$ is bounded, whereas the two latter ones assume
that the approximation of the gradient of $\mathsf{g}$ is unbiased,
that is
\begin{equation}
\label{e:unbiased}
(\forall n\in\NN)\quad\EC{u_n}{\XX_n}=\nabla \mathsf{g}(x_n).
\end{equation}

\section{Stochastic primal-dual splitting}\label{sec:PD}

The subdifferential 
\begin{equation}
\label{e:subdiff}
\partial\mathsf{f}\colon\mathsf{x}\mapsto
\menge{\mathsf{u}\in\HH}{(\forall\mathsf{y}\in\HH)\;\;
\scal{\mathsf{y}-\mathsf{x}}{\mathsf{u}}+\mathsf{f}(\mathsf{x})
\leq\mathsf{f}(\mathsf{y})}
\end{equation}
of a function  $\mathsf{f}\in \Gamma_0(\HH)$ is an example of a
maximally monotone operator \cite{Livre1}. Forward-backward
splitting has been developed in the more general framework 
of solving monotone inclusions \cite{Opti04,Livre1}. 
This powerful framework makes it possible 
to design efficient primal-dual strategies for optimization
problems; see for instance 
\cite{Opti14,Komo15} and the references therein. More 
precisely, we are interested in the following optimization problem
\cite[Section~4]{Svva12}.

\begin{problem}
\label{prob:3}
Let $\mathsf{f}\in \Gamma_0(\HH)$, let $\mu^{-1}\in\RPP$, let
$\mathsf{h}\colon\HH\to\RR$ be convex and differentiable with a 
$\mu^{-1}$-Lipschitz-continuous gradient,
and let $q$ be a strictly positive integer. For every 
$k\in\{1,\ldots,q\}$, let $\GG_k$ be a finite-dimensional 
Hilbert space, let 
$\mathsf{g}_k\in\Gamma_0(\GG_k)$,  
and let $\mathsf{L}_{k}\colon\HH\to\GG_k$ be linear. 
Let $\boldsymbol{\GG}=\GG_1\oplus\cdots\oplus\GG_q$ be the
direct Hilbert sum of $\GG_1,\ldots,\GG_q$, and suppose that
there exists $\overline{\mathsf{x}}\in\HH$ such that 
\begin{equation}\label{e:qualif}
0\in\partial\mathsf{f}(\overline{\mathsf{x}})
+\sum_{k=1}^q\mathsf{L}_{k}^*\partial\mathsf{g}_k(\mathsf{L}_{k}\overline{\mathsf{x}})
+\nabla\mathsf{h}(\overline{\mathsf{x}}).
\end{equation}
Let $\mathsf{F}$ be the set of solutions to the problem
\begin{equation}
\label{e:primopt}
\minimize{\mathsf{x}\in\HH}
{\mathsf{f}(\mathsf{x})+
\sum_{k=1}^q \mathsf{g}_k
(\mathsf{L}_{k}\mathsf{x})}+\mathsf{h}(\mathsf{x})
\end{equation}
and let $\boldsymbol{\mathsf{F}}^*$
be the set of  solutions to the dual problem 
\begin{equation} \label{e:dualopt}
\minimize{\boldsymbol{\mathsf{v}}\in\boldsymbol{\GG}}
{(\mathsf{f}^*\infconv \mathsf{h}^*) \bigg(-\Sum_{k=1}^q
\mathsf{L}_{k}^*\mathsf{v}_{k}\bigg)+\sum_{k=1}^q
\mathsf{g}_k^*(\mathsf{v}_k)},
\end{equation}
where $\infconv$ denotes the infimal convolution operation 
and $\boldsymbol{\mathsf{v}}=
(\mathsf{v}_1,\ldots,\mathsf{v}_q)$ designates a generic 
point in $\boldsymbol{\GG}$.
The objective is to find a point in 
$\mathsf{F}\times\boldsymbol{\mathsf{F}}^*$.
\end{problem}

We are interested in the case when only stochastic
approximations of the gradients of $\mathsf{h}$ and approximations
of the function $\mathsf{f}$ are available to solve
Problem~\ref{prob:3}. The following algorithm, which can be 
viewed as a stochastic extension of those of
\cite{Bang13,Cham11,Esse10,Heyu12,Icip14,Opti14,Cond13},
will be the focus of our investigation.

\begin{algorithm}
\label{algo:7}
Let $\rho \in \RPP$, let $(\mathsf{f}_n)_{n\in\NN}$ be a sequence
of functions in $\Gamma_0(\HH)$,
let $(\lambda_n)_{n\in\NN}$ be a sequence in $\left]0,1\right]$ 
such that $\sum_{n\in\NN}\lambda_n=\pinf$, and, for every 
$k\in\{1,\ldots,q\}$, let $\sigma_k\in\RPP$.
Let $x_0$, $(u_n)_{n\in\NN}$, and
$(b_n)_{n\in\NN}$ be random variables in $L^2(\Omega,\FF,\PP;\HH)$,
and let $\boldsymbol{v}_0$
and $(\boldsymbol{c}_n)_{n\in\NN}$ be random variables in 
$L^2(\Omega,\FF,\PP;\GGG)$.
Iterate
\begin{equation}
\label{e:PDcoordopt1}
\begin{array}{l}
\text{for}\;n=0,1,\ldots\\
\left\lfloor
\begin{array}{l}
\displaystyle y_{n} =
\prox_{\rho\mathsf{f}_n}
\left(x_{n}-\rho\bigg(\sum_{k=1}^q {\mathsf{L}^*_{k} v_{k,n}}
+u_n\bigg)\right)+b_{n}\\
x_{n+1}=x_{n}+\lambda_n (y_{n}-x_{n})\\
\text{for}\;k=1,\ldots,q\\
\left\lfloor
\begin{array}{l}
\displaystyle
w_{k,n}=\prox_{\sigma_k\mathsf{g}_{k}^*}
\big(v_{k,n}+\sigma_k\mathsf{L}_{k}
(2y_{n}-x_{n})\big)+c_{k,n}\\
v_{k,n+1}=v_{k,n}+\lambda_n (w_{k,n}-v_{k,n}).
\end{array}
\right.
\end{array}
\right.\\
\end{array}
\end{equation}
\end{algorithm}

One of main benefits of the proposed algorithm is that
it allows us to solve jointly the primal problem \eqref{e:primopt}
and the dual one \eqref{e:dualopt} in a fully decomposed fashion,
where each function and linear operator is activated individually.
In particular, it does not require any inversion of some
linear operator related to the operators $(\mathsf{L}_k)_{1\le k
\le q}$ arising in the original problem.  The convergence of the
algorithm is guaranteed by the following result which follows from
\cite[Proposition~5.3]{Pafa1}.

\begin{proposition}
\label{p:3}
Consider the setting of Problem~\ref{prob:3}, let 
$\mathscr{X}=(\boldsymbol{\XX}_n)_{n\in\NN}$ be a sequence of
sub-sigma-algebras of $\FF$, and let $(x_n)_{n\in\NN}$ and 
$(\boldsymbol{v}_n)_{n\in \NN}$ be sequences generated by 
Algorithm~\ref{algo:7}. 
Suppose that the following are satisfied:
\begin{enumerate}[label=\rm(\alph*)]
\item \label{a:p30}
$(\forall n\!\in\!\NN)$ 
$\sigma(x_{n'},\boldsymbol{v}_{n'})_{0\leq n'\leq n}\subset
\XXX_n\subset\boldsymbol{\XX}_{n+1}$.
\item 
\label{a:p3i}
$\sum_{n\in\NN}\lambda_n\sqrt{\EC{\| b_n\|^2}
{\boldsymbol{\XX}_n}}<\pinf$ and\newline
$\sum_{n\in\NN}\lambda_n\sqrt{\EC{\| \boldsymbol{c}_n\|^2}
{\boldsymbol{\XX}_n}}<\pinf$.
\item 
\label{a:p3ii}
$\sum_{n\in\NN}\sqrt{\lambda_n}\|\EC{u_n}{\boldsymbol{\XX}_n}-
\nabla\mathsf{h}(x_{n}) \|<\pinf$.
\item 
\label{a:p3iv}
There exists a summable sequence $(\tau_n)_{n\in\NN}$ in $\RP$ such 
that, for every $\mathsf{x}\in \mathsf{F}$, there exists 
$\big(\zeta_n(\mathsf{x})\big)_{n\in\NN}\in
\ell^\infty_+({\mathscr{X}})$ such that
$\big(\lambda_n\zeta_n(\mathsf{x})\big)_{n\in\NN}\in
\ell_+^{1/2}({\mathscr{X}})$ and
\begin{multline}
\label{e:boundtauetazetaPD}
(\forall n\in\NN)\quad 
\EC{\|u_n-\EC{u_n}{\boldsymbol{\XX}_n}\|^2}{\boldsymbol{\XX}_n}\\
\leq\tau_n \|\nabla\mathsf{h}(x_n)-\nabla\mathsf{h}
(\mathsf{x})\|^2
+\zeta_n(\mathsf{x}).
\end{multline}
\item 
\label{a:p3v}
There exist sequences $(\alpha_{n})_{n\in \NN}$ and 
$(\beta_{n})_{n\in \NN}$ in $\RP$ such that 
$\sum_{n\in \NN} \sqrt{\lambda_n} \alpha_{n} < \pinf$, 
$\sum_{n\in\NN} \lambda_n \beta_{n}$ $<\pinf$, and
\begin{multline}
(\forall n\in\NN)(\forall \mathsf{x}\in\HH)\\
\|\prox_{\rho\mathsf{f}_n}\mathsf{x}-
\prox_{\rho \mathsf{f}}\mathsf{x}\| 
\leq\alpha_{n} \|\mathsf{x}\|+\beta_{n}.
\end{multline}
\item  
\label{a:p3vi}
$\left(\rho^{-1}-\sum_{k=1}^q \sigma_k 
\|\mathsf{L}_{k} \|^2\right)\mu > 1/2$.
\end{enumerate}
Then, for some $\mathsf{F}$-valued random
variable $x$ and some $\boldsymbol{\mathsf{F}}^*$-valued random 
variable $\boldsymbol{v}$, $(x_n)_{n\in\NN}$ converges almost
surely to $x$ and
$(\boldsymbol{v}_n)_{n\in\NN}$ converges almost surely
to $\boldsymbol{v}$. 
\end{proposition}

\section{Application to online signal recovery}
\label{se:appli}
We consider the recovery of a signal 
$\overline{\mathsf{x}}\in\HH=\RR^N$ from the observation
model
\begin{equation}
(\forall n \in \NN)\quad z_{n} = K_{n} \overline{\mathsf{x}}+e_{n},
\end{equation}
where $K_n$ is a $\RR^{M\times N}$-valued random matrix and 
$e_n$ is a $\RR^M$-valued random noise vector. The objective is to
recover $\overline{\mathsf{x}}$ from $(K_n,z_n)_{n\in \NN}$, which
is assumed to be an identically distributed sequence. Such
recovery problems have been addressed in \cite{Comb89}. 
In this context, we propose to solve the primal 
problem \eqref{e:primopt} with $q=1$ and
\begin{equation}
(\forall \mathsf{x} \in \RR^N)\qquad \mathsf{h}(\mathsf{x})=\frac12
\EE\|K_0\mathsf{x}-z_0 \|^2,
\end{equation}
while functions $\mathsf{f}$ and $\mathsf{g}_1\circ \mathsf{L}_1$
are used to promote prior information on the target solution.
Since the statistics of the sequence $(K_n,z_n)_{n\in \NN}$ are not
assumed to be known a priori and have to be learnt online, at
iteration $n\in \NN$, we employ the empirical
estimate 
\begin{equation}
u_n=\frac{1}{m_{n+1}}\sum_{n'=0}^{m_{n+1}-1} 
K_{n'}^\top (K_{n'} x_n-z_{n'})
\end{equation}
of $\nabla \mathsf{h}(x_n)$.
The following statement, which can be deduced from
\cite[Section~5.2]{Pafa1}, illustrates the applicability of 
the results of Section~\ref{sec:PD}.

\begin{proposition}
Consider the setting of Problem~\ref{prob:3} and
Algorithm~\ref{algo:7}, where $\mathsf{f}_n\equiv \mathsf{f}$,
$b_n\equiv 0$, and $\boldsymbol{c}_n \equiv 0$. Let
$(m_n)_{n\in\NN}$ be a strictly increasing sequence in $\NN$ such
that $m_n=O(n^{1+\delta})$ with $\delta\in\RPP$, and let
\begin{equation}
(\forall n\in\NN)\quad\boldsymbol{\XX}_n=
\sigma\big(x_0,\boldsymbol{v}_0,(K_{n'},e_{n'})_{0\leq n'<m_n}).
\end{equation}
Suppose that the following are satisfied:
\begin{enumerate}[label=\rm(\alph*)]
\item 
The domain of $\mathsf{f}$ is bounded.
\item 
$(K_n,e_n)_{n\in \NN}$, is an independent and 
identically distributed (i.i.d.) sequence such that
$\EE\|K_0\|^4<\pinf$ and\linebreak 
$\EE\|e_0\|^4<\pinf$.
\item  
$\lambda_n=O(n^{-\kappa})$, where 
$\kappa\in\left]1-\delta,1\right]\cap [0,1]$.
\end{enumerate}
Then Assumptions~\ref{a:p30}-\ref{a:p3v} in Proposition~\ref{p:3}
hold.
\end{proposition}

Based on this result, we apply Algorithm~\ref{algo:7} 
to a practical scenario in which a grayscale image of
size $256\times 256$ with pixel values in $[0,255]$ is
degraded by a stochastic blur. The stochastic operator 
corresponds to a uniform i.i.d. subsampling of a uniform $5\times
5$ blur, performed in the discrete Fourier domain. More precisely,
the value of the frequency response at each frequency bin is
kept with probability $0.3$ or it is set to zero. In addition, the
image is corrupted by an additive zero-mean white Gaussian 
noise with standard deviation equal to $5$. 
The average signal-to-noise ratio (SNR) is initially equal to 
$3.4$~dB. 
\begin{figure}[htb]
\begin{center}
\includegraphics[width=7.0cm]{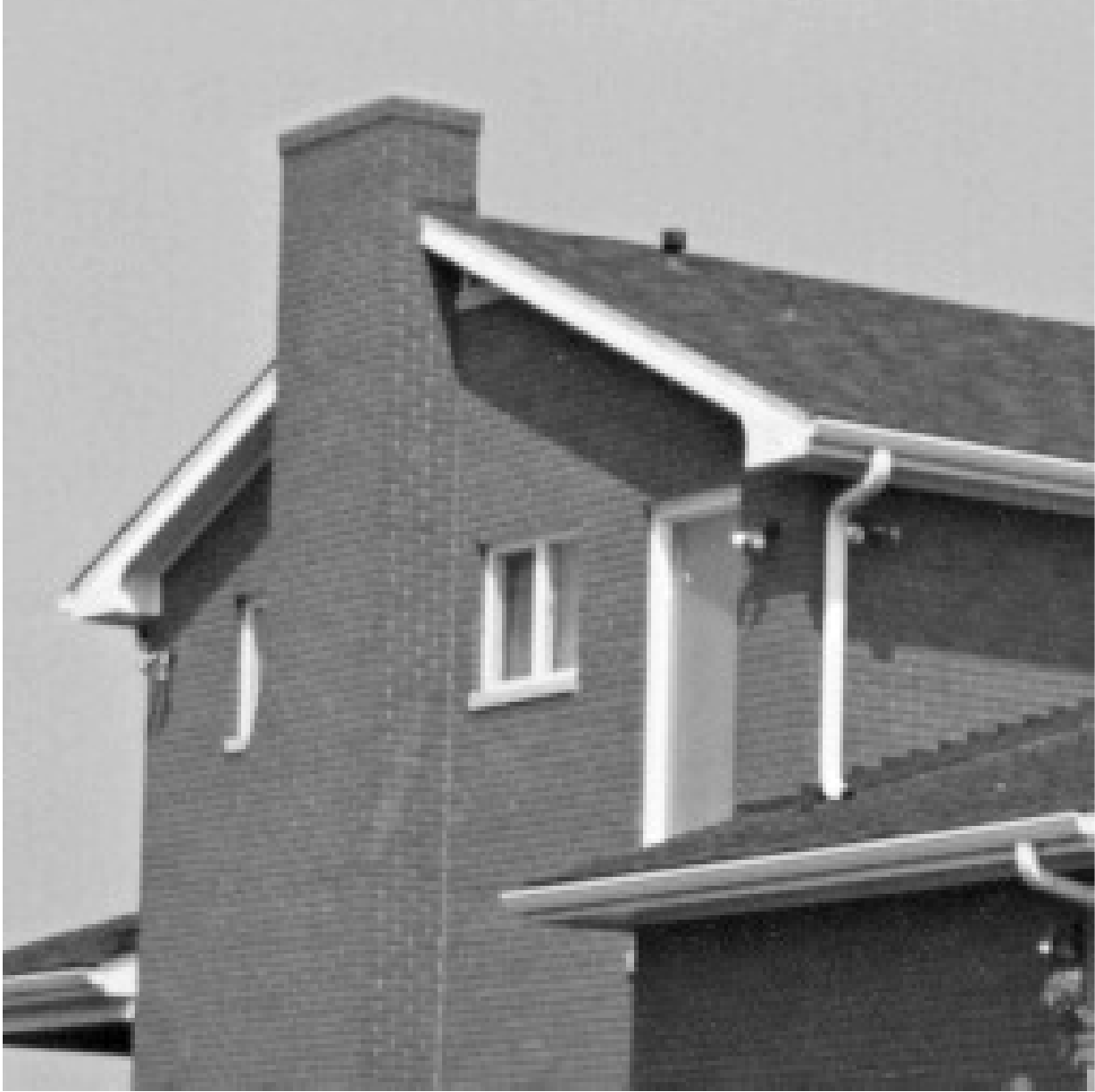}\\
\includegraphics[width=7.0cm]{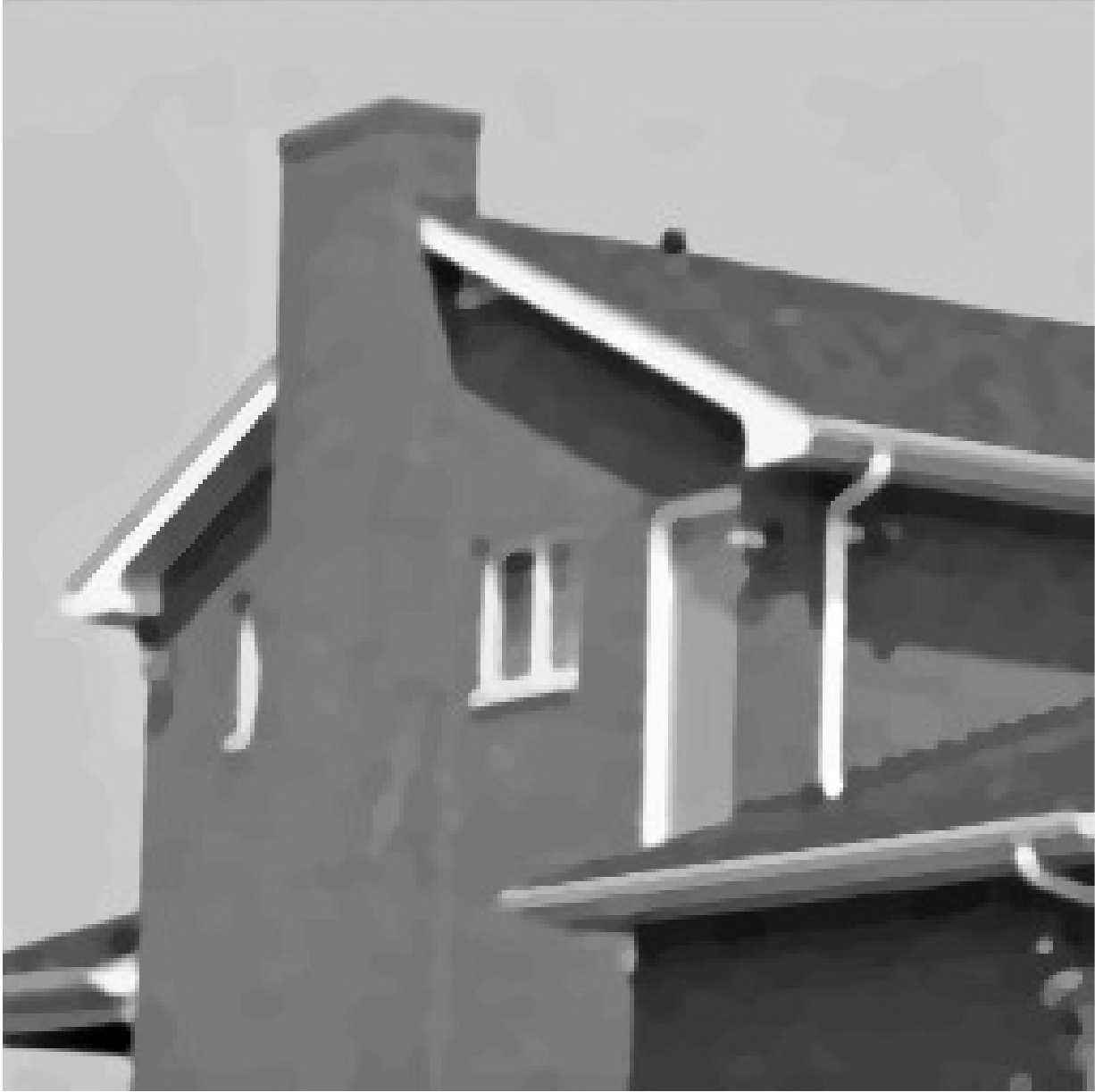}\\
\end{center}
\caption{Original image $\overline{x}$ (top), 
restored image (bottom).}
\label{fig:1}
\end{figure}

\begin{figure}[htb]
\begin{center}
\includegraphics[width=7.0cm]{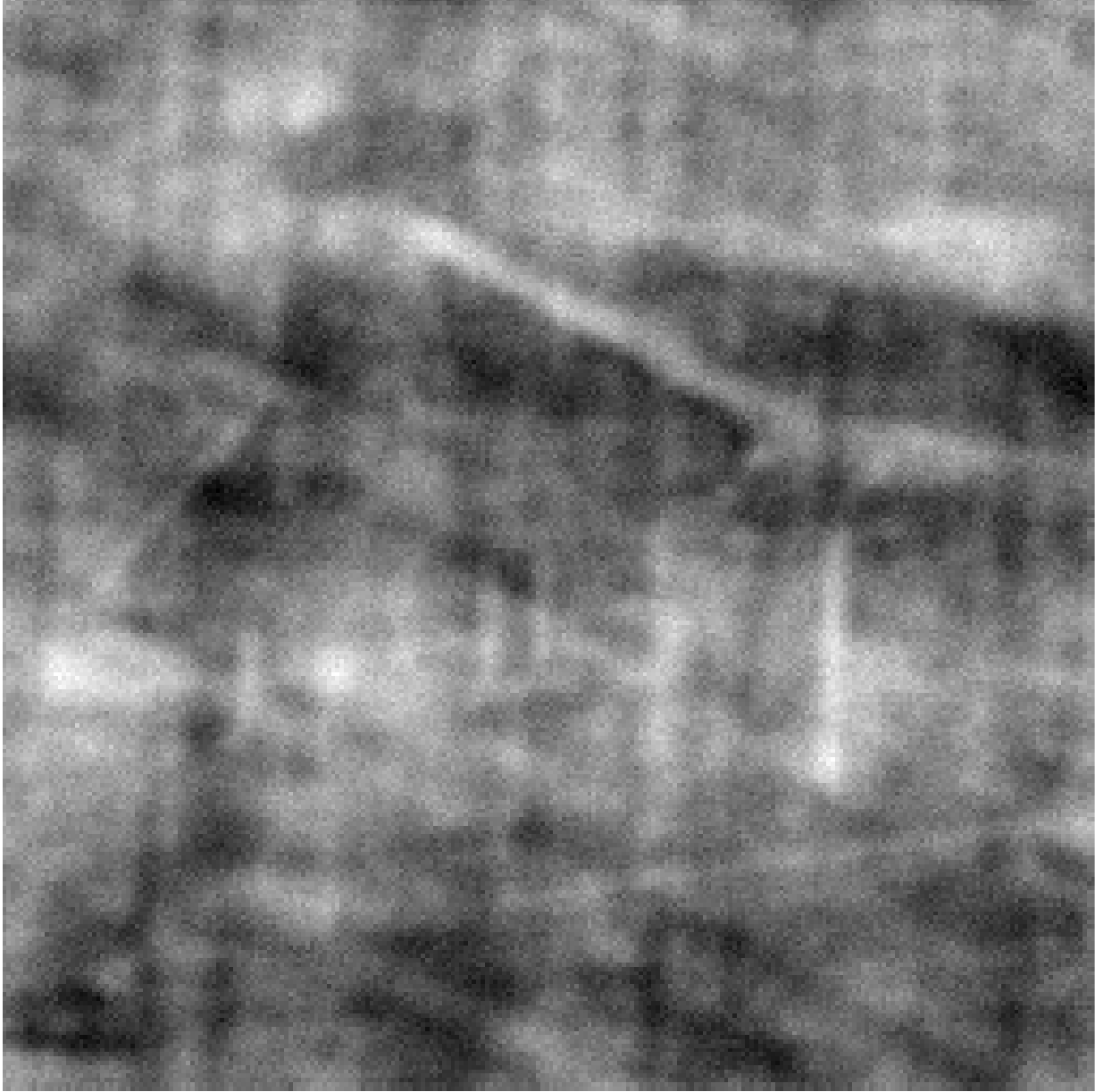}\\
\includegraphics[width=7.0cm]{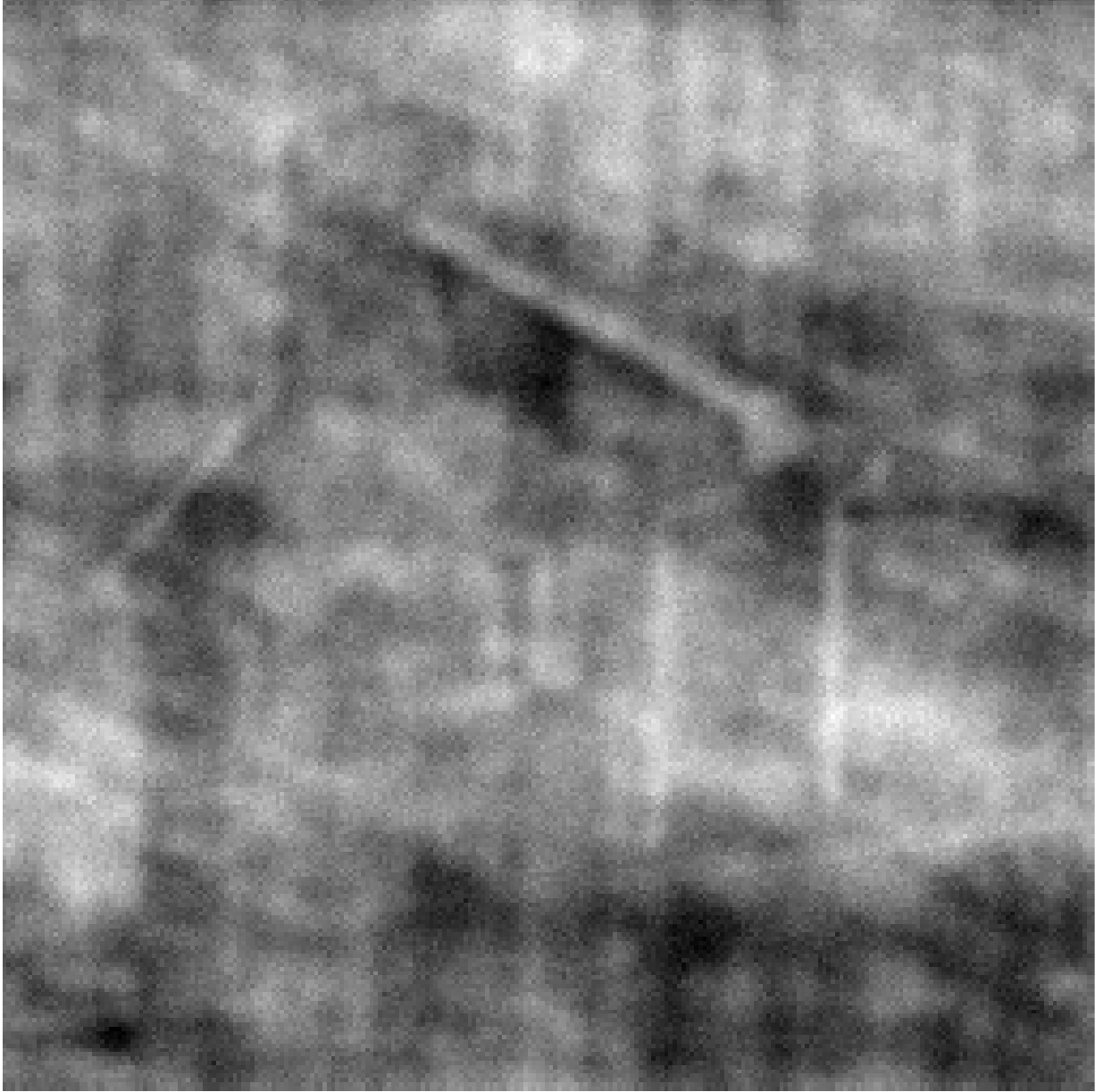}\\
\end{center}
\caption{
Degraded image 1 (SNR = $0.14$~dB) (top), and degraded image 2
(SNR = $12.0$~dB)  (bottom).}
\label{fig:1bis}
\end{figure}

In our restoration approach, the function $\mathsf{f}$ is
the indicator function of the set
$[0,255]^N$, while $\mathsf{g}_1\circ
\mathsf{L}_1$ is a classical isotropic total variation regularizer,
where $\mathsf{L}_1$ is the concatenation of the horizontal and
vertical discrete gradient operators. 
Figs.~\ref{fig:1}--\ref{fig:1bis} displays the
original image, the restored image, 
as well as two realizations of the
degraded images. The SNR for the restored image is equal to
$28.1$~dB.

\begin{figure}[htb]
\centering
\includegraphics[width=8cm]{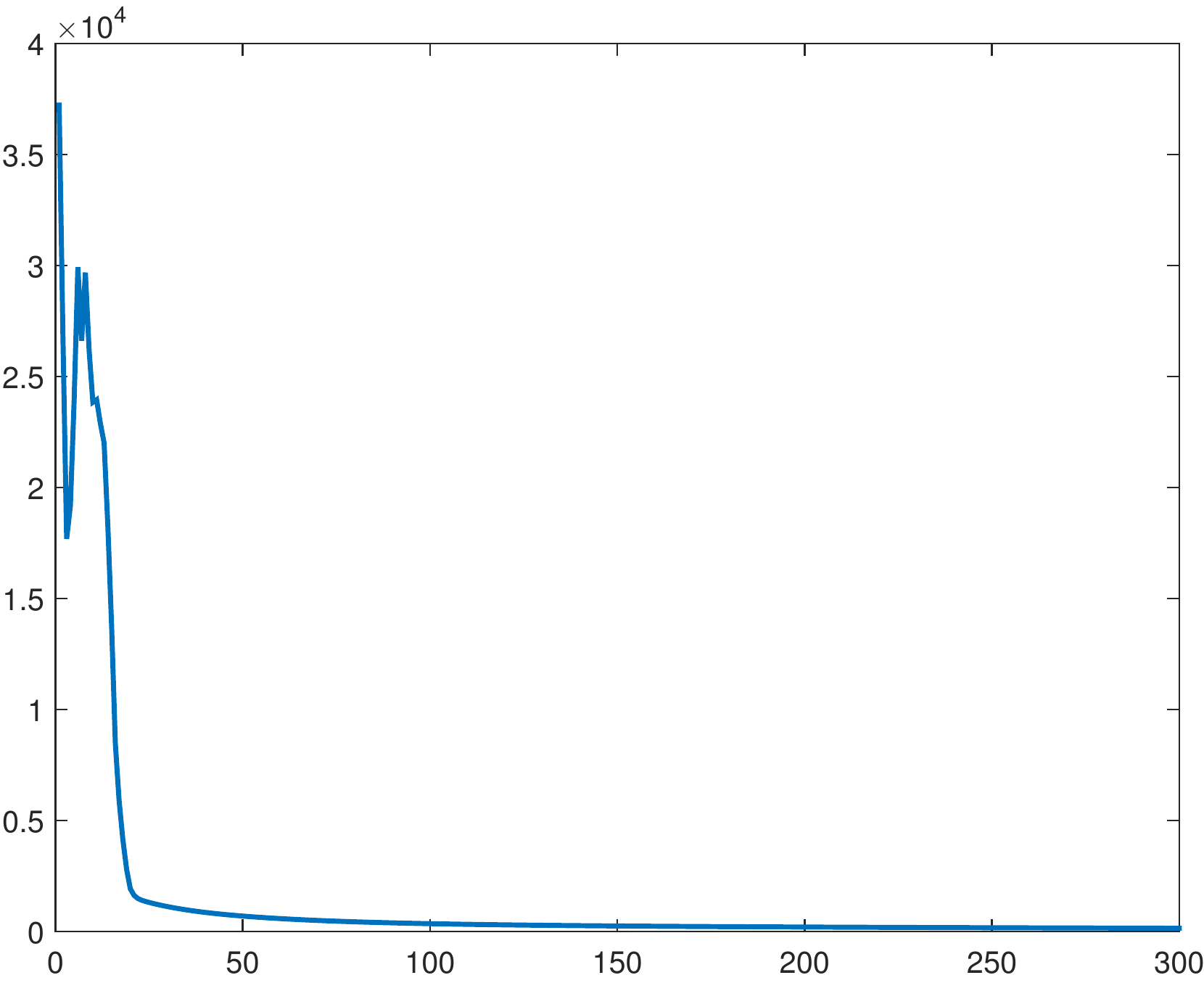}
\caption{\small 
$\|x_n-x_\infty\|$ versus the iteration number $n$.}
\label{fig:2}
\end{figure}

Fig.~\ref{fig:2} shows the convergence behavior
of the algorithm. In these experiments, we have chosen
\begin{equation}
(\forall n \in \NN)\quad \begin{cases}
m_n = n^{1.1}\\
\lambda_n = (1+(n/500)^{0.95})^{-1}.
\end{cases}
\end{equation}

\section{Conclusion}\label{se:conclu}
We have proposed two stochastic proximal splitting
algorithms for solving nonsmooth convex optimization
problems. These methods require only approximations of the
functions used in the formulation of the optimization problem,
which is of the utmost importance for solving online signal
processing problems. The almost sure convergence of these algorithms
has been established. The stochastic
version of the primal-dual algorithm that we have investigated has
been evaluated in an online image restoration problem in which the
data are blurred by a stochastic point spread function and 
corrupted with noise.

\bibliographystyle{IEEEbib}

\end{document}